\documentclass[a4paper,12pt,reqno]{amsart}
\usepackage[text={130mm,200mm}]{geometry}

\theoremstyle{plain}

\newtheorem*{theorem*}{Theorem}

\newtheorem*{corollary*}{Corollary}
\newtheorem{lemma}{Lemma}
\newtheorem*{lemma*}{Lemma}
\newtheorem{proposition}{Proposition}
\newtheorem*{proposition*}{Proposition}

\newtheorem*{conjecture*}{Conjecture}
\theoremstyle{definition}

\newtheorem*{definition*}{Definition}
\theoremstyle{remark}

\newtheorem*{remark*}{Remark}

\begin{document}

\title[Modeling rational numbers by Cantor series ]{Modeling rational numbers by Cantor series}
\author{Symon Serbenyuk}
\address{
  45~Shchukina St. \\
  Vinnytsia \\
  21012 \\
  Ukraine}
\email{simon6@ukr.net}



\begin{abstract}
  In the present article, modeling certain rational numbers, that are represented  in terms of Cantor series,  are described. The statements on relations between digits in the representations  of  rational numbers by Cantor series (for the case of an arbitrary sequence $(q_k)$) are proved.
\end{abstract}

\maketitle

\maketitle



\section{Introduction}

Let $Q\equiv (q_k)$ be a fixed sequence of positive integers, $q_k>1$,  $\Theta_k$ be a sequence of the sets $\Theta_k\equiv\{0,1,\dots ,q_k-1\}$, and $\varepsilon_k\in\Theta_k$.

The Cantor series expansion 
\begin{equation}
\label{eq: Cantor series}
\frac{\varepsilon_1}{q_1}+\frac{\varepsilon_2}{q_1q_2}+\dots +\frac{\varepsilon_k}{q_1q_2\dots q_k}+\dots
\end{equation}
of $x\in [0,1]$,   first studied by G. Cantor in \cite{Cantor1}. It is easy to see that the Cantor series expansion is the  b-ary expansion
$$
\frac{\alpha_1}{b}+\frac{\alpha_2}{b^2}+\dots+\frac{\alpha_n}{b^n}+\dots
$$
of numbers  from the closed interval $[0,1]$ whenever the condition $q_k=b$ holds for all positive integers $k$. Here $b$ is a fixed positive integer, $b>1$, and $\alpha_n\in\{0,1,\dots , b-1\}$.

By $x=\Delta^Q _{\varepsilon_1\varepsilon_2\ldots\varepsilon_k\ldots}$  denote a number $x\in [0,1]$ represented by series \eqref{eq: Cantor series}. This notation is called \emph{the representation of $x$ by Cantor series \eqref{eq: Cantor series}.}

We note that certain numbers from $[0,1]$ have two different representations by Cantor series \eqref{eq: Cantor series}, i.e., 
$$
\Delta^Q _{\varepsilon_1\varepsilon_2\ldots\varepsilon_{m-1}\varepsilon_m000\ldots}=\Delta^Q _{\varepsilon_1\varepsilon_2\ldots\varepsilon_{m-1}[\varepsilon_m-1][q_{m+1}-1][q_{m+2}-1]\ldots}=\sum^{m} _{i=1}{\frac{\varepsilon_i}{q_1q_2\dots q_i}}.
$$
Such numbers are called \emph{$Q$-rational}. The other numbers in $[0,1]$ are called \emph{$Q$-irrational}.

Let $c_1,c_2,\dots, c_m$ be an
ordered tuple of integers such that $c_i\in\{0,1,\dots, q_i-~1\}$ for $i=\overline{1,m}$. 

\emph{A cylinder $\Delta^Q _{c_1c_2...c_m}$ of rank $m$ with base $c_1c_2\ldots c_m$} is a set of the form
$$
\Delta^Q _{c_1c_2...c_m}\equiv\{x: x=\Delta^Q _{c_1c_2...c_m\varepsilon_{m+1}\varepsilon_{m+2}\ldots\varepsilon_{m+k}\ldots}\}.
$$
That is any cylinder $\Delta^Q _{c_1c_2...c_m}$ is a closed interval of the form
$$
\left[\Delta^Q _{c_1c_2...c_m000}, \Delta^Q _{c_1c_2...c_m[q_{m+1}][q_{m+2}][q_{m+3}]...}\right].
$$

Define \emph{the shift operator $\sigma$ of expansion \eqref{eq: Cantor series}} by the rule
$$
\sigma(x)=\sigma\left(\Delta^Q _{\varepsilon_1\varepsilon_2\ldots\varepsilon_k\ldots}\right)=\sum^{\infty} _{k=2}{\frac{\varepsilon_k}{q_2q_3\dots q_k}}=q_1\Delta^{Q} _{0\varepsilon_2\ldots\varepsilon_k\ldots}.
$$

It is easy to see that 
\begin{equation*}
\label{eq: Cantor series 2}
\begin{split}
\sigma^n(x) &=\sigma^n\left(\Delta^Q _{\varepsilon_1\varepsilon_2\ldots\varepsilon_k\ldots}\right)\\
& =\sum^{\infty} _{k=n+1}{\frac{\varepsilon_k}{q_{n+1}q_{n+2}\dots q_k}}=q_1\dots q_n\Delta^{Q} _{\underbrace{0\ldots 0}_{n}\varepsilon_{n+1}\varepsilon_{n+2}\ldots}.
\end{split}
\end{equation*}

Therefore, 
\begin{equation}
\label{eq: Cantor series 3}
x=\sum^{n} _{i=1}{\frac{\varepsilon_i}{q_1q_2\dots q_i}}+\frac{1}{q_1q_2\dots q_n}\sigma^n(x).
\end{equation}

Note that, in the paper \cite{Serbenyuk2017}, the notion of the shift operator of an alternating Cantor series is studied in detail.

In 1864, G.~Cantor introduced the problem on representaions of rational numbers by series \eqref{eq: Cantor series} (see \cite{Cantor1}). More information about this problem were described in \cite{SA17}. 

In the present article, we consider the case of positive Cantor series. In the next articles, the cases of alternating and sign-variable Cantor series and certain applications of used techniques will be  consider by the author of the present article.

\section{Representations of certain rational numbers}

Using the main statements from \cite{Cantor1, S13, Ser2017} (the paper \cite{S13} is  \cite{Ser2017}  translated into English), we get the following. 
 
\begin{proposition}
A rational number $x=\frac{p}{r}$, where $p<r$ and (p,r)=1,  has a finite expansion by positive Cantor series  whenever there exists a number $k_0$ such that the condition  $q_1q_2\cdots q_{k_0}\equiv 0 \pmod{r}$ holds. 
\end{proposition}

\begin{proposition}
There exist certain sequences $(q_k)$ such that all rational numbers represented in terms of  corresponding Cantor series have  finite expansions. 
\end{proposition}
For example, these representations are following:
$$
x=\Delta^{(2k)} _{\varepsilon_1\varepsilon_2...\varepsilon_k...}\equiv\sum^{\infty} _{k=1}{\frac{\varepsilon_k}{2\cdot 4\cdot 8\cdot \ldots \cdot 2k}},~\text{where}~\varepsilon_k\in\{0,1,\dots,2k-1\};
$$
$$
x=\Delta^{(k+1)!} _{\varepsilon_1\varepsilon_2...\varepsilon_k...}\equiv\sum^{\infty} _{k=1}{\frac{\varepsilon_k}{2\cdot 3\cdot 4\cdot \ldots \cdot (k+1)}},~\text{where}~\varepsilon_k\in\{0,1,\dots,k\}.
$$

\begin{proposition}
Suppose that the sequence $(q_k)$ is periodic. Then a number $x$ is rational if and only if the representation $\Delta^Q _{\varepsilon_1\varepsilon_2\ldots\varepsilon_k\ldots}$ of $x$ is periodic (i.e., the sequence $(\varepsilon_k)$ is periodic). 
\end{proposition}

\begin{proposition}
The following is true:
$$
\frac{1}{w}=\sum^{\infty} _{k=1}{\frac{\varepsilon^{'} _{k}}{q_1q_2\cdots q_k}},
$$
where  $w\varepsilon^{'} _{k}=q_k-1$ holds for all positive integers $k$, $w$ is a certain positive integer.
\end{proposition}
It follows from the condition $1=\Delta^Q _{[q_1-1][q_2-1]...[q_k-1]...}$.

\begin{proposition}
Let $n_0$ be a fixed positive integer number, $q_0=\min_{n>n_0}{q_n}$, and $\varepsilon_0$ be a numerator of the fraction 
$\frac{\varepsilon_{n_0+k}}{q_1q_2...q_{n_0}q_{n_0+1}...q_{n_0+k}}$ in expansion \eqref{eq: Cantor series} of $x$ providing that $q_{n_0+k}=q_0$. Then $\sigma^n(x)=const$ for all $n\ge n_0$ if and only if the condition $\frac{q_n-1}{q_0-1}\varepsilon_0=\varepsilon_n\in\mathbb Z_0$ holds for any $n>n_0$.
\end{proposition}

\begin{proposition}
A number $x$ represented by expansion \eqref{eq: Cantor series} is rational if and only if there exists a subsequence $(n_k)$ of positive integers  such that for all $k=1,2,\dots ,$ the following conditions are true:
\begin{itemize}
\item
$$
\frac{\lambda_k}{\mu_k}=
\frac{\varepsilon_{n_k+1}q_{n_k+2}\dots q_{n_{k+1}}+\varepsilon_{n_k+2}q_{n_k+3}\dots q_{n_{k+1}}+\dots +\varepsilon_{n_{k+1}-1}q_{n_{k+1}}+\varepsilon_{n_{k+1}}}{q_{n_k+1}q_{n_k+2}\dots q_{n_{k+1}}-1}=const;
$$
\item $\lambda_k=\frac{\mu_k}{\mu}\lambda$, where $\mu=\min_{k\in\mathbb N}{\mu_k}$ and $\lambda$ is a number in the numerator of the fraction whose denominator equals $(\mu_1+1)(\mu_2+~1)\dots (\mu+~1)$  from sum \eqref{eq: Cantor series 6}.
\end{itemize}
Here
$$
x=\sum^{\infty} _{k=1}{\frac{\varepsilon_k}{q_1q_2\dots q_k}}=\sum^{n_1} _{j=1}{\frac{\varepsilon_j}{q_1q_2\dots q_j}}+\frac{1}{q_1q_2\dots q_{n_1}}x^{'},
$$
$$
x^{'}=\sum^{\infty} _{k=1}{\frac{\varepsilon_{n_k+1}q_{n_k+2}q_{n_k+3}\dots q_{n_{k+1}}+\varepsilon_{n_k+2}q_{n_k+3}\dots q_{n_{k+1}}+\dots+\varepsilon_{n_{k+1}-1}q_{n_{k+1}}+\varepsilon_{n_{k+1}}}{(q_{n_1+1}\dots q_{n_2})(q_{n_2+1}\dots q_{n_3})\dots (q_{n_k+1}\dots q_{n_{k+1}})}}
$$
\begin{equation}
\label{eq: Cantor series 6}
=\sum^{\infty} _{k=1}{\frac{\lambda_k}{(\mu_1+1)\dots (\mu_k+1)}}.
\end{equation}
\end{proposition}
For example, 
$$
x=\Delta^{(2k+1)} _{123...k...}\equiv\frac{1}{3}+\frac{2}{3\cdot 5}+\frac{3}{3\cdot 5\cdot 7}+\dots +\frac{k}{3\cdot 5\cdot \ldots \cdot (2k+1)}+\dots =\frac{1}{2}.
$$
Here $x=\sigma^n(x)=\frac{1}{2}$ for all $n$, where $n=~0,1,2, \dots$.

The last-mentioned sum is useful for modeling rational numbers of the type $\frac{A}{B}$, where $B\equiv 0 \pmod{2}$. In particular,
$$
\frac{1}{6}=\Delta^{(2k+1)} _{0234...k...},~\frac{5}{6}=\Delta^{(2k+1)} _{2234...k...}.
$$

\section{The main results}

Let $\frac{p}{r}$ be a fixed number, where $(p,r)=1$, $p<r$, and $p\in\mathbb N, r\in\mathbb N$. Here $\mathbb N$ is the set of all positive integers. Then
$$
\frac{p}{r}=\sum^{\infty} _{n=1}{\frac{\delta_n}{q_1q_2\cdots q_n}}.
$$
\begin{remark*}
We use the denotation $\delta_n$ for the known $n$th digit in the representation of a number $x$ by a Cantor series and $\varepsilon_n$ for the  non-known $n$th digit. 

In addition, since $ x\in \Delta^Q _{c_1c_2...c_m}$
but $\Delta^Q _{c_1c_2...c_m[q_{m+1}-1][q_{m+2}]...}=\Delta^Q _{c_1c_2...c_{m-1}[c_m+1]000...}$,
we assume that
$$
\Delta^Q _{c_1c_2...c_{m-1}[c_m]000...}\le x<\Delta^Q _{c_1c_2...c_{m-1}[c_m+1]000...}.
$$
\end{remark*}
 
It is easy to see that
$$
\frac{p}{r}\in\Delta^Q _{\varepsilon_1}=\left[\Delta^Q _{\varepsilon_1000...}, \Delta^Q _{\varepsilon_1[q_2-1][q_3-1]...}\right]=\left[\frac{\varepsilon_1}{q_1},\frac{\varepsilon_1+1}{q_1}\right].
$$
That is 
$$
\frac{\varepsilon_1}{q_1}\le \frac{p}{r}< \frac{\varepsilon_1+1}{q_1},
$$
$$
\varepsilon_1r\le pq_1<r(\varepsilon_1+1),
$$
$$
\varepsilon_1\le \frac{pq_1}{r}<\varepsilon_1+1.
$$
So,
$$
\varepsilon_1=\left[\frac{p}{r}q_1\right]\equiv \delta_1,
$$
where $[x]$ is the integer part of $x$. 

Now we get 
$$
\frac{p}{r}\in \Delta^Q _{\delta_1\varepsilon_2}=\left[\frac{q_2\delta_1+\varepsilon_2}{q_1q_2},\frac{q_2\delta_1+\varepsilon_2+1}{q_1q_2}\right].
$$
Whence,
$$
\frac{q_2\delta_1+\varepsilon_2}{q_1q_2}\le \frac{p}{r}< \frac{q_2\delta_1+\varepsilon_2+1}{q_1q_2},
$$
$$
\varepsilon_2\le \frac{pq_1q_2-rq_2\delta_1}{r}<\varepsilon_2+1.
$$
So,
$$
\varepsilon_2=\left[\frac{pq_1q_2-rq_2\delta_1}{r}\right]\equiv\delta_2. 
$$

In the third step, we have
$$
\frac{p}{r}\in \Delta^Q _{\delta_1\delta_2\varepsilon_3}=\left[\frac{\delta_1q_2q_3+\delta_2q_3+\varepsilon_3}{q_1q_2q_3},\frac{\delta_1q_2q_3+\delta_2q_3+\varepsilon_3+1}{q_1q_2q_3}\right]
$$
and
$$
\varepsilon_3=\left[\frac{pq_1q_2q_3-r(\delta_1q_2q_3+\delta_2q_3)}{r}\right]\equiv\delta_3.
$$

In the $n$th step, we obtain
$$
\frac{p}{r}\in \Delta^Q _{\delta_1\delta_2...\delta_{n-1}\varepsilon_n}=\left[\sum^{n-1} _{i=1}{\frac{\delta_i}{q_1q_2\cdots q_i}}+\frac{\varepsilon_n}{q_1q_2\cdots q_n},\sum^{n-1} _{i=1}{\frac{\delta_i}{q_1q_2\cdots q_i}}+\frac{\varepsilon_n+1}{q_1q_2\cdots q_n}\right]
$$
$$
=\left[\frac{\delta_1q_2q_3\cdots q_n+\delta_2q_3q_4\cdots q_n+\dots+\delta_{n-1}q_n+\varepsilon_n}{q_1q_2\cdots q_n},\frac{\delta_1q_2q_3\cdots q_n+\dots+\delta_{n-1}q_n+\varepsilon_n+1}{q_1q_2\cdots q_n}\right].
$$

Let $\varsigma_n$ denote the sum $\delta_1q_2q_3\cdots q_n+\delta_2q_3q_4\cdots q_n+\dots+\delta_{n-1}q_n$. Then
$$
\frac{\varsigma_n+\varepsilon_n}{q_1q_2\cdots q_n}\le \frac{p}{r}<\frac{\varsigma_n+\varepsilon_n+1}{q_1q_2\cdots q_n},
$$
$$
\varepsilon_n\le\frac{pq_1q_2\cdots q_n-r\varsigma_n}{r}<\varepsilon_n+1.
$$
Denoting by $\Delta_n=pq_1q_2\cdots q_n-r\varsigma_n$, we get
$$
\varepsilon_n=\left[\frac{\Delta_n}{r}\right]\equiv \delta_n.
$$
So, the following statement is true.

\begin{lemma}
Let $x\in (0,1)$ be a rational number represented by series~\eqref{eq: Cantor series}.
If $x=\frac{p}{q}=\Delta^Q _{\delta_1\delta_2...\delta_n...}$, then the equality
\begin{equation}
\label{eq: delta 1}
\delta_n=\left[\frac{\Delta_n}{r}\right]
\end{equation}
 holds for all $n\in\mathbb N$,
where
$$
\Delta_n=pq_1q_2\cdots q_n-r(\delta_1q_2q_3\cdots q_n+\delta_2q_3q_4\cdots q_n+\dots+\delta_{n-1}q_n).
$$
\end{lemma}

Also, for $n\ge 2$ the condition $\varsigma_n=\varsigma_{n-1}q_n+\delta_{n-1}q_n$ holds and
$$
\Delta_n=q_n(\Delta_{n-1}-r\delta_{n-1}).
$$

\begin{lemma}
Let $x\in (0,1)$ be a rational number represented by series~\eqref{eq: Cantor series}.
If $x=\frac{p}{q}=\Delta^Q _{\delta_1\delta_2...\delta_n...}$, then the equality
$$
\delta_n=\left[\frac{q_n(\Delta_{n-1}-r\delta_{n-1})}{r}\right]
$$
 holds for all $1<n\in\mathbb N$, where $\Delta_1=pq_1$ and $\delta_1=\left[\frac{\Delta_1}{r}\right]$.

\end{lemma}

Suppose that the following sequence of conditions is true:
$$
\delta_1=\left[\frac{p}{r}q_1\right],~~~\delta_2=\left[\frac{pq_1q_2-rq_2\delta_1}{r}\right],~~~\delta_3=\left[\frac{pq_1q_2q_3-r(\delta_1q_2q_3+\delta_2q_3)}{r}\right], \dots 
$$
$$
\dots , \delta_n=\left[\frac{pq_1q_2\cdots q_n}{r}-(\delta_1q_2q_3\cdots q_n+\delta_2q_3q_4\cdots q_n+\dots+\delta_{n-1}q_n)\right]=\left[\frac{\Delta_n}{r}\right], \dots .
$$

It follows from equality \eqref{eq: Cantor series 3} that
$$
x=\frac{p}{r}=\frac{\delta_1q_2q_3\cdots q_n+\delta_2q_3q_4\cdots q_n+\dots+\delta_{n-1}q_n+\delta_n}{q_1q_2\cdots q_n}+\frac{\sigma^n\left(\frac{p}{r}\right)}{q_1q_2\cdots q_n},
$$
$$
\sigma^n\left(\frac{p}{r}\right)=\frac{\Delta_n-r\delta_n}{r}=\frac{\Delta_n}{r}-\delta_n
$$
and
$$
\delta_n=\frac{\Delta_n}{r}-\sigma^n\left(\frac{p}{r}\right).
$$
From the last-mentioned relationship and relationship \eqref{eq: delta 1} it follows that
$$
0\le \delta_n=\left[\frac{\Delta_n}{r}\right]=\frac{\Delta_n}{r}-\sigma^n\left(\frac{p}{r}\right)=\left[\frac{\Delta_n}{r}-\sigma^n\left(\frac{p}{r}\right)\right].
$$
That is 
$$
\sigma^n\left(\frac{p}{r}\right)=\left\{\frac{\Delta_n}{r}\right\},
$$
where $\{a\}$ is the fractional part of $a$.

\begin{remark*}
Clearly, 
$$
0 \le \sigma^n\left(x\right)\le 1
$$
for an arbitrary $x\in [0,1]$. However for any Q-rational number $x=\Delta^Q _{\delta_1\delta_2...\delta_{n-1}\delta_n000...}=\Delta^Q _{\delta_1\delta_2...\delta_{n-1}[\delta_n-1][q_{n+1}-1][q_{n+2}-1][q_{n+3}-1]...}$ the following conditions hold:
$$
\sigma^n\left(x\right)=\sigma^n\left(\Delta^Q _{\delta_1\delta_2...\delta_{n-1}\delta_n000...}\right)=\sum^{\infty} _{k=n+1}{\frac{0}{q_{n+1}q_{n+2}\cdots q_{k}}}=0,
$$
$$
\sigma^n\left(x\right)=\sigma^n\left(\Delta^Q _{\delta_1\delta_2...\delta_{n-1}[\delta_n-1][q_{n+1}-1][q_{n+2}-1][q_{n+3}-1]...}\right)=\sum^{\infty} _{k=n+1}{\frac{q_k-1}{q_{n+1}q_{n+2}\cdots q_{k}}}=1.
$$
Since the condition  $\sigma^n\left(x\right)=1$ holds only for the last-mentioned representations of Q-rational numbers $x$, we can use only the first representation of Q-rational numbers $x$ and for these numbers the condition $\sigma^n\left(x\right)=0$ holds.
\end{remark*}

In addition, note that
$$
\varsigma_n=\frac{pq_1q_2\cdots q_n-\Delta_n}{r}.
$$
Whence for an arbitrary $n\in\mathbb N$
$$
x=\sum^{n} _{k=1}{\frac{\delta_k}{q_1q_2\cdots q_k}}+\frac{\sigma^n\left(x\right)}{q_1q_2\cdots q_n}
$$
$$
=\frac{\delta_1q_2q_3\cdots q_n+\delta_2q_3q_4\cdots q_n+\dots +\delta_{n-1}q_n+\delta_n}{q_1q_2\cdots q_n}+\frac{\sigma^n\left(x\right)}{q_1q_2\cdots q_n}
$$
$$
=\frac{\varsigma_n+\delta_n}{q_1q_2\cdots q_n}+\frac{\left\{\frac{\Delta_n}{r}\right\}}{q_1q_2\cdots q_n}=\frac{\frac{pq_1q_2\cdots q_n-\Delta_n}{r}+\delta_n}{q_1q_2\cdots q_n}+\frac{\frac{\Delta_n}{r}-\left[\frac{\Delta_n}{r}\right]}{q_1q_2\cdots q_n}
$$
$$
=\frac{pq_1q_2\cdots q_n-\Delta_n+r\delta_n}{rq_1q_2\cdots q_n}+\frac{\frac{\Delta_n}{r}-\delta_n}{q_1q_2\cdots q_n}=\frac{p}{r}.
$$

So, we have the following statement.
\begin{theorem*}
A number $x=\Delta^Q _{\delta_1\delta_2...\delta_n...} \in (0,1)$ is a rational number $\frac{p}{r}$, where $p,r\in\mathbb N, (p,r)=1$, and $p<r$, if and only if the condition 
$$
\delta_n=\left[\frac{q_n(\Delta_{n-1}-r\delta_{n-1})}{r}\right]
$$
 holds for all $1<n\in\mathbb N$, where $\Delta_1=pq_1$, $\delta_1=\left[\frac{\Delta_1}{r}\right]$, and 
 $[a]$ is the integer part of $a$.

\end{theorem*}

Let us consider certain examples. Suppose
$$
x=\Delta^{(2n+1)} _{\varepsilon_1\varepsilon_2...\varepsilon_n...}=\sum^{\infty} _{n=1}{\frac{\varepsilon_n}{3\cdot 5\cdot 7 \cdot \ldots \cdot (2n+1)}}.
$$
This representation is complicated for modeling rational numbers with an even denominator:
$$
\frac{1}{4}=\Delta^{(2n+1)} _{035229[11]4...}
$$
$$
\frac{3}{8}=\Delta^{(2n+1)} _{104341967...}.
$$

Investigations of the author of the present article about representations of rational numbers by Cantor series can be useful for solving ``P vs NP Problem" (this problem described in http://www.claymath.org/millennium-problems/p-vs-np-problem,  www.claymath.org/sites/default/files/pvsnp.pdf). The next articles of the author of the present article will be devoted to such investigations.


\begin{thebibliography}{9}

\bibitem{Cantor1} G.~Cantor, Ueber die einfachen Zahlensysteme, 
\emph{Z. Math. Phys. }   {\bf 14} (1869), 121--128.

\bibitem{Serbenyuk2017}
Symon Serbenyuk. Representation of real numbers by the alternating Cantor series, \emph{ Integers} {\bf  17}(2017), Paper No. A15, 27 pp.

\bibitem{S13}  S. Serbenyuk.  Cantor series and rational numbers, available at  https://arxiv.org/pdf/1702.00471.pdf

\bibitem{SA17}  S. Serbenyuk. Cantor series expansions of rational numbers, arXiv:1706.03124 or available at https://www.researchgate.net/publication/317099134

\bibitem{Ser2017}
S.~ Serbenyuk. Rational numbers in terms of positive Cantor series,  \emph{Bull. Taras Shevchenko Natl. Univ. Kyiv Math. Mech.} \textbf{36} (2017), no.~2,  11--15 (in Ukrainian)









\end{thebibliography}
\end{document}